\documentclass[10pt]{article}

\usepackage{amsmath}

\usepackage{graphicx}        % standard LaTeX graphics tool

\usepackage{labelfig}

\newtheorem{theorem}{Theorem} 
\newtheorem{remark}{Remark} 
\newtheorem{definition}{Definition}
\newtheorem{proposition}{Proposition} 
\newtheorem{lemma}{Lemma}

\begin{document}

\baselineskip=4.6mm

\makeatletter

\newcommand{\E}{\mathrm{e}\kern0.2pt}
\newcommand{\D}{\mathrm{d}\kern0.2pt}
\newcommand{\RR}{\mathrm{I\kern-0.20emR}}

\def\bottomfraction{0.9}

\title{{\bf A comparison theorem for super- and subsolutions of $\mathbf{\nabla^2 u
+ f (u) = 0}$ and its application to water waves with vorticity}}

\author{Vladimir Kozlov$^a$, Nikolay Kuznetsov$^b$}

\date{}

\maketitle

\vspace{-8mm}

\begin{center}
$^a${\it Department of Mathematics, Link\"oping University, S--581 83 Link\"oping,
Sweden \\ $^b$ Laboratory for Mathematical Modelling of Wave Phenomena, \\ Institute
for Problems in Mechanical Engineering, Russian Academy of Sciences, \\ V.O.,
Bol'shoy pr. 61, St. Petersburg 199178, RF} \\

\vspace{2mm}

E-mail: vladimir.kozlov@liu.se ; nikolay.g.kuznetsov@gmail.com
\end{center}

\begin{abstract}
A comparison theorem is proved for a pair of solutions that satisfy in a weak sense
opposite differential inequalities with nonlinearity of the form $f (u)$ with $f$
belonging to the class $L^p_{loc}$. The solutions are assumed to have non-vanishing
gradients in the domain, where the inequalities are considered. The comparison
theorem is applied to the problem describing steady, periodic water waves with
vorticity in the case of arbitrary free-surface profiles including overhanging ones.
Bounds for these profiles as well as streamfunctions and admissible values of the
total head are obtained. \vspace{2mm}

\noindent {\bf Keywords:} Comparison theorem, nonlinear differential inequality,
partial hodograph transform in $n$ dimensions, periodic steady water waves with
vorticity, streamfunction
\end{abstract}

\section{Introduction}

\setcounter{equation}{0}

In their remarkable article \cite{GNN}, Gidas, Ni and Nirenberg investigated various
properties that solutions (in particular, positive solutions) of several nonlinear
equations have in different (bounded as well as unbounded) domains in $\RR^n$, $n
\geq 1$. For this purpose several forms of the maximum principle were employed along
with some other methods. The authors emphasised that their techniques could be
applicable in physical situations other than those considered in the paper. During
the decades past since the publication of \cite{GNN}, this prediction proved
correct. The most spectacular results obtained in the paper deal with the equation
\begin{equation}
\nabla^2 u + f (u) = 0 , \quad \nabla u = (u_{x_1},\dots,u_{x_n}) \ \mbox{and} \
u_{x_i} = \partial_i u = \partial u / \partial x_i , \label{eq}
\end{equation}
in which $f$ is a $C^1$-function.

The two-dimensional version of \eqref{eq} describes, in particular, periodic steady
water waves with vorticity in which case $f$ is the given vorticity distribution. If
the depth of water is finite, the domain is a quadrangle bounded by three straight
segments\,--\,two of them that are opposite to each other are equal in view of
periodicity\,--\,and a curve that is opposite the third segment and corresponds to
the smallest period of wave propagating on the free surface; of course, one can
consider a strip with periodic upper boundary and horizontal bottom as the water
domain. (The relevant free-boundary problem is derived from Euler's equations, for
example, in \cite{CS}.)

Instead of equation \eqref{eq}, the present paper deals with the inequality
\begin{equation} 
\nabla^2 u + f (u) \leq 0 \quad \mbox{in a domain} \ \Omega \subset \RR^n ,
\ n \geq 2 , \label{in}
\end{equation}
and its opposite which are understood in a weak sense. In the case of \eqref{in},
this means that the integral inequality
\begin{equation} 
\int_X \nabla u \cdot \nabla v \, \D x \geq \int_X f (u) v \, \D x \label{inX}
\end{equation}
is valid for every non-negative $v \in C_0^1 (X)$, where $X$ is any subdomain of
$\Omega$. Moreover, no smoothness and even continuity is required from $f$, and the
aim is to prove the following comparison theorem for a pair of functions that
satisfy the inequalities.

\begin{theorem}
Let $f \in L^p_{loc} (\RR)$ with $p > n$, and let $u_1, u_2 \in C^{1} (\Omega)$
have non-va\-nishing gradients in $\Omega$ and satisfy in the weak sense the
inequalities
\begin{equation} 
\nabla^2 u_1 + f (u_1) \geq 0 \quad and \quad \nabla^2 u_2 + f (u_2) \leq 0 \quad in
\ \Omega , \label{ineq}
\end{equation}
respectively. If $u_1 \leq u_2$ in $\Omega$ and these functions are equal at some
point $x^0 \in \Omega$, then $u_1$ and $u_2$ coincide throughout $\Omega$.
\end{theorem}

\begin{remark}
{\rm If $f \in L^p (\Omega)$ and $u \in C^{1} (\Omega)$ is such that $\nabla u \neq
0$ throughout $\Omega$, then $f (u (x))$ is a measurable function in $\Omega$;
moreover, this superposition belongs to $L^p_{loc} (\Omega)$.

It should be mentioned that Theorem~1 is not true without the assumption that the
gradients of $u_1$ and $u_2$ are non-vanishing. Indeed, even for H\"older continuous
$f$ (the weaker condition $f \in L^p (\Omega)$, $p > n$, is imposed in Theorem~1),
this follows from the example on p.~220 in \cite{GNN}.

Let $\Omega = \RR^n$ and $u_1$ be equal to zero identically. If $p > 2$, then $u_2$
equal to $(1 - |x|^2)^p$ when $|x| \leq 1$ and to zero for $|x| > 1$ belongs to $C^2
(\RR^n)$. It is straightforward to check that \eqref{eq} holds for $u_2$ with
\[ f (u) = - 2 p (p-2) u^{1-2/p} + 2 p (n + 2 p - 2) u^{1-1/p} , 
\]
which is H\"older continuous with the exponent $1 - 2/p$ and such that $f (0) = 0$.
Thus, all assumptions of Theorem~1 are fulfilled for $u_1$ and $u_2$ except for the
condition concerning their gradients; for both functions they vanish when $|x| \geq 1$.
Therefore, the conclusion of Theorem~1 is not true\,--\,these functions do not
coincide.}
\end{remark}

The proof of Theorem 1 is given in \S\,2; it is based on the so-called partial
hodograph transform in $n$ dimensions which allows us to use the weak Harnack type
inequality proved in \cite{T}. In \S\,3, we apply Theorem 1 to obtain bounds for
solutions of the free-boundary problem mentioned above; it describes steady,
periodic water waves with vorticity.

\section{Proof of Theorem 1}

First, the local version of the $n$-dimensional partial hodograph transform is
introduced. It is used in the proof of an auxiliary lemma required for proving
Theorem 1. Then a version of Hopf's lemma is discussed; the latter is applied in
considerations of \S\,3.

\subsection{The partial hodograph transform in $\mathbf{n}$ dimensions}

Being defined locally, it generalises the transform introduced by Dubreil-Jacotin
\cite{DJ} in her studies of water waves with vorticity in the two-dimensional case.
Moreover, the transform proposed here extends that considered in \cite{DJ} to the
case of $n > 2$ dimensions. Therefore, it might be of interest for applications
other than water waves with vorticity.

Let $\Omega \subset \RR^n$ be a domain. If $u$ is a function whose gradient
does not vanish throughout this domain, then at any point $x^0 \in \Omega$ the
coordinate system can be chosen so that $u_{x_n} (x^0) > 0$. This allows us to
introduce the following transform in a neighbourhood of $x^0$. We put
\[ q = (q_1,\dots,q_{n-1}) \ \mbox{with} \ q_k = x_k, \ k=1,\dots,n-1, \
\mbox{and} \ p = u (x) ,
\]
and take these as new independent variables. Furthermore, instead of $u (x)$
satisfying, say inequality \eqref{in}, we consider $h (q, p) = x_n$ as the unknown.
Then we have
\[ \frac{\partial h}{\partial x_k} = h_{q_k} + h_p \frac{\partial u}{\partial x_k} = 0
\ \mbox{for} \ k=1,\dots,n-1 \ \mbox{and} \ \frac{\partial h}{\partial x_n} = h_p
u_{x_n} = 1 , 
\]
and so
\[ h_{q_k} = - \frac{u_{x_k}}{u_{x_n}} \ \mbox{for} \ k=1,\dots,n-1 \ \mbox{and}
\ h_p = \frac{1}{u_{x_n}} > 0 \, .
\]
In view of the equalities
\[ u_{x_k} = - \frac{h_{q_k}}{h_p} , \ k=1,\dots,n-1 , \ \ \mbox{and} \
u_{x_n} = \frac{1}{h_p} \, ,
\]
the weak formulation \eqref{inX} of inequality \eqref{in} for $u$ takes the
following form in terms of $h$:
\begin{equation} 
\int_Q \left[ - \nabla_q h \cdot \nabla_q w + \frac{1 + |\nabla_q h|^2}{h_p} w_p
\right] \D q \D p \leq \int_Q f (p) w h_p \, \D q \D p \, . \label{inQ}
\end{equation}
Here $Q$\,--\,a neighbourhood of the point $(q^0, p^0)$\,--\,is the image of $X$
which is the neighbourhood of $x^0$ that corresponds to $(q^0, p^0)$ and $w (q, p)$
stands instead of $v (x (q, p))$. It is also taken into account that 
\[ \D x = h_p \,
\D q \D p \quad \mbox{and} \quad \nabla_q h = (h_{q_1}, \dots,h_{q_{n-1}}) .
\]
Like \eqref{inX}, the last inequality must hold for every non-negative $w \in C_0^1
(Q)$.

In the case of smooth $h$, a consequence of \eqref{inQ} is the differential
inequality
\[ (L h) (q, p) \geq f(p) h_p (q, p) , \quad \mbox{where} \ L h =
\partial_{q_k}h_{q_k} - \partial_p \frac{1 + |\nabla_q h|^2}{h_p} \, .
\]
It should be noted that the right-hand side of this inequality depending on $f$ is
linear in $h$, whereas nonlinearity is present in the operator $L$ on the left-hand
side. It will be clear from what follows that this is the advantage resulting from
the introduced form of the hodograph transform. Moreover, this differs from what we
have in \eqref{in}, where the differential operator is linear and nonlinearity is
involved through the superposition operator $f (u)$.

Furthermore, let $h_1, h_2 \in C^{1} (Q)$, then we have
\[ L h_1 - L h_2 = \partial_{q_k} (h_1 - h_2)_{q_k} - \partial_p \int_0^1 \partial_t
\frac{1 + |\nabla_q h^{(t)}|^2}{h^{(t)}_p} \, \D t ,
\]
where $h^{(t)} = t h_1 + (1-t) h_2$. This difference can be written as an operator
in divergent form with continuous coefficients. Indeed, let $W = h_1 - h_2$, then $L
h_1 - L h_2 = \mathcal{L} W$, where
\[ \mathcal{L} W = \partial_{q_k} W_{q_k} - 2 \, \partial_p \! \left[ W_{q_k} \int_0^1
\frac{h^{(t)}_{q_k}}{h^{(t)}_p} \, \D t \right] + \partial_p \! \left[ W_p \int_0^1
\frac{1 + |\nabla_q h^{(t)}|^2}{h^{(t) 2}_p} \, \D t \right] .
\]
Moreover, the inequalities
\[ \sum_{k=1}^{n-1} \left[ \int_0^1 \frac{h^{(t)}_{q_k}}{h^{(t)}_p} \, \D t \right]^2
\leq \sum_{k=1}^{n-1} \int_0^1 \frac{h^{(t) 2}_{q_k}}{h^{(t) 2}_p} \, \D t <
\int_0^1 \frac{1 + |\nabla_q h^{(t)}|^2}{h^{(t) 2}_p} \, \D t 
\]
show that $\mathcal{L}$ is an elliptic operator.

\subsection{Auxiliary lemma}

This immediate corollary of Theorem~5.1 in \cite{T} is given for the reader's
convenience. It concerns an inequality for the linear elliptic operator $\mathcal{P}
(\partial) = \partial_j (a_{ij} \partial_i) + b_j \partial_j + c$. As usual the
ellipticity means that there exists $\lambda > 0$ such that
\[ \lambda^{-1} |\xi|^2 < a_{ij} \xi_i \xi_j < \lambda |\xi|^2 \quad \mbox{for all}
\ \xi = (\xi_1,\dots,\xi_n) \in \RR^n \setminus \{0\} .
\]
Here and below, the Einstein summation notation is used.

\begin{lemma}
Let $u \in C^{1} (\Omega)$ be a non-negative function satisfying the inequality
\begin{equation} 
\mathcal{P} (\partial) u \leq 0 \quad in \ \Omega \label{ineq'}
\end{equation}
in the weak sense; here all coefficients $a_{ij}$ are measurable in $\Omega$, $b_j
\in L^q (\Omega)$ and $c \in L^{q/2} (\Omega)$ for some $q > n$. If there exists
$x^0 \in \Omega$ such that $u (x^0) = 0$, then $u$ vanishes identically in $\Omega$.
\end{lemma}

\noindent {\it Proof.} Let $\rho > 0$ be such that $K_{3 \rho} (x^0) \subset
\Omega$; by $K_{\sigma} (x)$ the open cube centred at $x$ is denoted which has edges
equal to $\sigma$ and sides parallel to the coordinate axes. The lemma's assumptions
about $u$ and $\mathcal{P}$ yield the inequality
\[ \rho^{-n/\gamma} \| u \|_{L^\gamma (K_{2 \rho} (x^0))} \leq C \min_{x \in 
K_{\rho} (x^0)} u (x)
\]
for $u$ satisfying \eqref{ineq'}; here $C$ is a positive constant and $\gamma$ is an
arbitrary number from the interval $(1, n / (n-2))$. Since
\[ \min_{x \in K_{\rho} (x^0)} u (x) = u (x^0) = 0 \, ,
\]
there exists a neighbourhood of $x^0$, where $u$ vanishes identically. It is clear
that the maximal such neighbourhood is $\Omega$ because otherwise the same argument
can be applied to any boundary point of the maximal neighbourhood which is an
interior point of $\Omega$, thus leading to a contradiction.

\begin{figure}[t!]
\centering
\vspace{-2mm}%
\SetLabels
 \L (-0.01*0.94) $p$\\
 \L (0.95*-0.05) $x_n$\\
 \L (-0.02*0.485) $p^0$\\
 \L (0.77*0.83) $u_2 (x)$\\
 \L (0.86*0.57) $u_1 (x)$\\
 \L (0.2*-0.065) $h_2 (q, p^0)$\\
 \L (0.54*-0.065) $h_1 (q, p^0)$\\
\endSetLabels \leavevmode \strut\AffixLabels{\includegraphics[width=60mm]{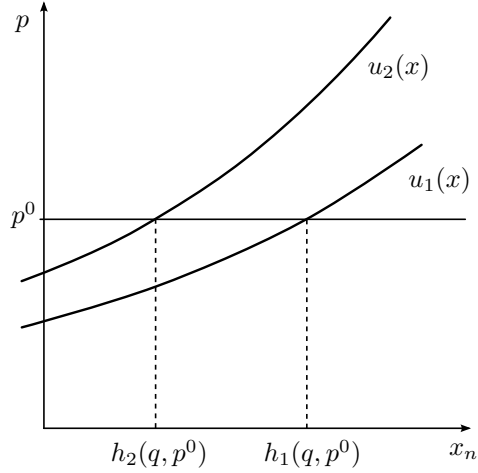}}
\vspace{2mm} \caption{A sketch of the partial hodograph transform with fixed $x_1 =
q_1, \, \dots, \, x_{n-1} = q_{n-1}$. It demonstrates the relationship between
inequalities for two pairs of corresponding functions. If $u_1 (x) < u_2 (x)$ near
$(x_1,\dots,x_{n-1},x_n^0)$, then $h_1(q, p) > h_2(q, p)$ in a neighbourhood of $(q,
p^0)$.} \vspace{-2mm}
\end{figure}

\begin{remark}
For a non-positive $u$ satisfying the inequality $\mathcal{P} (\partial) u \geq 0$
the assertion of Lemma 1 remains valid.
\end{remark}

\subsection{Proof of Theorem 1}

Let us apply the partial hodograph transform to $u_1$ and $u_2$ in a neighbourhood
of $x^0$. Then there exists a neighbourhood $Q$ of $(q^0, p^0)$\,--\,the image of
$x^0$\,--\,in which the inequalities
\begin{equation} 
(-1)^j \big[ (L h_j) (q, p) - f (p) \, \partial_p h_j (q, p) \big] \geq 0 , \quad
j=1,2, \label{jin}
\end{equation}
follow from \eqref{in}; here $h_1$ and $h_2$ are the functions corresponding to
$u_1$ and $u_2$, respectively. It is clear that $h_1 (q^0, p^0) = h_2 (q^0, p^0)$,
and $Q$ can be taken so that $h_1 \geq h_2$ in it because $u_1 \leq u_2$ in $\Omega$
(see Figure~1).

From \eqref{jin} one obtains that $W = h_1 - h_2$ (it is non-negative in $Q$)
satisfies the inequality
\[ (\mathcal{L} W) (q, p) - f (p) W (q, p) \leq 0 \quad \mbox{for} \ (q, p) \in Q .
\]
Then Lemma 1 yields that $W$ vanishes identically in $Q$, that is, $h_1 = h_2$
throughout $Q$, and so $u_1 = u_2$ in some neighbourhood of $x^0$. Thus the set,
say $E$, where $u_1$ coincides with $u_2$, is non-empty. It is clear that $E$ is
closed in $\Omega$, and so if $E \neq \Omega$, then there exists $x^* \in \partial E
\cap \Omega$. Applying the same considerations, we see that $x^*$ has a
neighbourhood belonging to $E$ which is a contradiction.

\subsection{Hopf's lemma}

In \S\,3 we will apply the following version of the well-known result whose proof
we failed to find in the literature.

\begin{lemma}
Let $x^0$ be a point on a part of $\partial \Omega$ belonging to the class
$C^{1,\alpha}$, $\alpha \in (0, 1)$, and let $u \in C^{1} (\bar X)$ satisfy
\eqref{ineq'} in the weak sense in $X$ which is the intersection of a neighbourhood
of $x^0$ with $\Omega$, whereas $\mathcal{P}$ is such that all $a_{ij} \in
C^{0,\alpha} (X)$, $b_j \in L^q (X)$ and $c \in L^{q/2} (X)$ for some $q > n$.

If $u$ is non-negative in $X$ and $u (x^0) = 0$, then either $u$ vanishes
identically in $X$ or $\partial_n u(x^0) < 0$, where $\partial_n$ denotes the normal
derivative on $\partial \Omega$ directed to the exterior of $\Omega$.
\end{lemma}

\noindent {\it Proof.} The proof is essentially the same as that of Theorem 1.1 in
\cite{Lis}, where the assumptions imposed on $b_j$ and $c$ (boundedness and
non-negativity of the last coefficient) are superfluous. Therefore, we restrict
ourselves to some necessary remarks and begin with a couple of minor notes. First,
it is sufficient to prove the assertion in the case when $\partial \Omega$ is flat
near $x_0$, to which the general case reduces by a change of variables. Second, in
view of Lemma~1, one has just to show that $\partial_n u(x^0) < 0$ is a consequence
of the inequality $u > 0$ in $X$.

The only amendment that needs more details concerns the proof of Lem\-ma~3.4 in
\cite{Lis}, in which $\nabla w_3$ should be estimated as follows (we keep the
notation used in \cite{Lis}):
\begin{eqnarray*}
|\nabla w_3 (x)| \leq C \|\nabla v\|_{L^\infty} \int_{\Omega_\rho} \frac{|b (y)| + c
(y) |x-y|}{|x-y|^{n-1}} \D y \\ \leq C \|\nabla v\|_{L^\infty} \big( \rho^{\alpha_1}
\|b\|_{L^q} + \rho^{\alpha_2} \|c\|_{L^{q/2}} \big) \ \, ,
\end{eqnarray*}
where $\alpha_1 = (q-n)/(q-1)$, $\alpha_2 = 2 (q-n)/(q-2)$.

\begin{remark}
{\rm For a non-positive $u$ satisfying the inequality $\mathcal{P} (\partial) u \geq
0$ the as\-sertion of Lemma 2 remains valid provided the conclusion that
\[ \partial_n u(x^0) < 0 \ \ \mbox{is changed to} \ \ \partial_n u(x^0) > 0. \]}
\end{remark}

\section{Application of Theorem 1 \\ to water waves with vorticity}

In this section, we consider the two-dimensional nonlinear problem of steady,
periodic waves in an open channel occupied by an inviscid, incompressible, heavy
fluid, say water. The water motion is assumed to be rotational which, according to
observations, is the type of motion commonly occurring in nature. A brief
characterization of results obtained earlier for this and other related problems is
given in our paper \cite{KK2}. Further details can be found in the survey article
\cite{WS} by Strauss; see also the recent papers \cite{CSV}, \cite{KKL1} and
\cite{KKL2}. Here, our aim is to apply Theorem~1 in order to generalise conditions
guaranteeing the validity of bounds for solutions to the problem that were obtained
in \cite{KK1}.

\subsection{Statement of the problem}

Let an open channel of uniform rectangular cross-section be bounded below by a
horizontal rigid bottom and let water occupying the channel be bounded above by a
free surface not touching the bottom. The surface tension is neglected and the
pressure is constant on the free surface. Since the water motion is supposed to be
two-dimensional and rotational and in view of the water incompressibility, we seek
the velocity field in the form $(\psi_y, -\psi_x)$, where the unknown function $\psi
(x,y)$ is referred to as the {\it streamfunction} (see, for example, \cite{LSh} for
details of this model). It is also supposed that the vorticity distribution $\omega$
(it is a function of $\psi$ as is explained in \cite{LSh}, \S\,1) is a prescribed
function belonging to $L^p_{loc} (\RR)$ with $p > 2$. This assumption is weaker
than that in our previous paper \cite{KK1}, where $\omega$ was assumed to be a
locally Lipschitz function.

Non-dimensional variables are chosen so that the constant volume rate of flow per
unit span and the constant acceleration due to gravity are scaled to unity in our
equations. In appropriate Cartesian coordinates $(x,y)$, the bottom coincides with
the $x$-axis and gravity acts in the negative $y$-direction. The frame of reference
is taken so that the velocity field is time-independent as well as the unknown
free-surface profile. The latter is assumed to be a simple $C^{1}$-curve, say
$\Gamma$, which is $\Lambda$-periodic along the $x$-axis for some $\Lambda > 0$, but
not necessarily representable as the graph of an $x$-dependent function. (It was
found numerically by Vanden-Broeck \cite{VB} that there are such overhanging
profiles bounding rotational flows with periodic waves, whereas Constantin, Strauss
and Varvaruca \cite{CSV} recently investigated them rigorously; the relevant figures
are presented in these papers.) Thus, the longitudinal section of the water domain
is the strip $D_\Gamma$ that lies between the $x$-axis and $\Gamma$, and $\psi
(x,y)$ is assumed to be a $\Lambda$-periodic function of $x$ in $D_\Gamma$.

Since the surface tension is neglected, the pair $(\psi, \Gamma)$ must be found from
the following free-boundary problem:
\begin{eqnarray}
&& \psi_{xx} + \psi_{yy} + \omega (\psi) = 0, \quad (x, y) \in D_\Gamma;
\label{eq:lapp} \\ && \psi (x,0) = 0, \quad x \in \RR; \label{eq:bcp} \\ && 
\psi (x, y) = 1, \quad (x, y) \in \Gamma; \label{eq:kcp} \\ && |\nabla \psi (x,
y)|^2 + 2 y = 3 r, \quad (x, y) \in \Gamma . \label{eq:bep}
\end{eqnarray}
Here $r$ is a constant considered as the given problem's parameter. Notice that the
boundary condition \eqref{eq:kcp} allows us to write relation \eqref{eq:bep}
(Bernoulli's equation) as follows:
\begin{equation} 
\left[ \partial_n \psi (x, y) \right]^2  + 2 y = 3 r, \quad (x, y) \in \Gamma \, .
\label{eq:ben}
\end{equation}
Here and below $\partial_n$ denotes the normal derivative on $\Gamma$; the normal $n
= (n_x, n_y)$ has unit length and points out of $D_\Gamma$. 

In this section, we keep the notation adopted in our previous papers, and so $\psi$
and $\omega$ stand in \eqref{eq:lapp} instead of $u$ and $f$, respectively, used in
\eqref{eq} and \eqref{in}. To give the precise definition how a solution of problem
\eqref{eq:lapp}--\eqref{eq:bep} is understood we need the following set $\Gamma_\psi
= \{ (x, y) \in D_\Gamma : \nabla \psi (x, y) = 0 \}$.

\begin{definition}
{\rm The pair $(\psi, \Gamma)$ is called a solution of problem
\eqref{eq:lapp}--\eqref{eq:bep} with the vorticity distribution $\omega \in
L^p_{loc} (\RR)$, $p > 2$, provided the following conditions are fulfilled for some
$\Lambda > 0$:

\noindent $\bullet$ $\Gamma$ is a simple, $\Lambda$-periodic along the $x$-axis
$C^{1}$-curve;

\noindent $\bullet$ $\psi (x,y)$ is a $\Lambda$-periodic function of $x$ belonging
to $C^{1} (\overline {D_\Gamma})$;

\noindent $\bullet$ the boundary conditions \eqref{eq:bcp}--\eqref{eq:bep} are
fulfilled pointwise;

\noindent $\bullet$ the two-dimensional measure of $\Gamma_\psi$ is equal to zero
and $D_\Gamma \setminus \Gamma_\psi$ is a domain;

\noindent $\bullet$ for all $v \in C_0^\infty (D_\Gamma \setminus \Gamma_\psi)$ the
following identity holds:
\[ \int_{D_\Gamma} \nabla \psi \cdot \nabla v \, \D x \D y = \int_{D_\Gamma} 
\omega (\psi) v \, \D x \D y .
\]}
\end{definition}

The last condition means that $\psi$ is a weak solution of \eqref{eq:lapp} in
$D_\Gamma \setminus \Gamma_\psi$.

\subsection{Auxiliary one-dimensional problems}

The results presented in this section were obtained in \cite{KK3} under the
assumption that $\omega$ is a Lipschitz function. Since they are essential for our
considerations, what follows is a digest of these results valid under the assumption
that $\omega \in L^1_{loc} (\RR)$.

First, let $s > 0$, then by $U (y; s)$ we denote a strictly monotonic solution of
the following Cauchy problem:
\[ U'' + \omega (U) = 0 , \ y \in \RR ; \ \ U (0; s) = 0 , \ \ U' (0; s) = s ;
\] 
here and below $'$ stands for $\D / \D y$. It is straightforward to obtain the
implicit formula
\begin{equation}
y = \int_0^U \frac{\D \tau}{\sqrt{s^2 - 2 \Omega (\tau)}} \, , \quad \Omega (\tau) =
\int_0^\tau \omega (t) \, \D t , \label{eq:Uim}
\end{equation}
that defines $U$ on the maximal interval of monotonicity $(y_- (s), y_+ (s))$,
where
\[ y_\pm (s) = \int_0^{\tau_\pm (s)} \frac{\D \tau}{\sqrt{s^2 - 2 \, \Omega (\tau)}}
\, ,
\]
and the definition of $\tau_\pm (s)$ is as follows. By $\tau_+ (s)$ and $\tau_- (s)$
we denote the least positive and the largest negative root, respectively, of the
equation $2 \, \Omega (\tau) = s^2$. If this equation has no positive (negative)
root, we put 
\[ \tau_+ (s) = +\infty \quad (\tau_- (s) = -\infty \ \mbox{respectively}) . \]

Second, we consider the problem
\begin{equation}
u'' + \omega (u) = 0 \ \ \mbox{on} \ (0, h) , \quad u (0) = 0 , \ \ u (h) = 1 ,
\label{eq:u}
\end{equation}
in the class of monotonic functions. It is clear that formula \eqref{eq:Uim} gives a
solution of problem \eqref{eq:u} on the interval $(0, h (s))$, where
\begin{equation}
h (s) = \int_0^1 \frac{\D \tau}{\sqrt{s^2 - 2 \, \Omega (\tau)}} \quad \mbox{and} \
s > s_0 = \sqrt{\, 2 \max_{\tau \in [0, 1]} \Omega (\tau)} \geq 0 .
\label{eq:d}
\end{equation}
Moreover, all monotonic solutions of problem \eqref{eq:u} have the form
\eqref{eq:Uim} on the interval $(0, h)$. This remains valid for $s = s_0$ with
\[ h = h_0 = \int_0^1 \frac{\D \tau}{\sqrt{s^2_0 - 2 \, \Omega (\tau)}} < \infty ,
\quad \mbox{that is}, \ h_0 = \lim_{s \to s_0} h (s) .
\]
It is clear that $h (s)$ decreases strictly monotonically from $h_0$ and asymptotes
zero as $s \to \infty$.

Furthermore, the pair $(u, \Gamma)$ with
\[ u (y) = U (y; s) \quad \mbox{and} \quad \Gamma = \{ (x, y) : x \in \RR , 
y = h (s) \}
\]
is a solution of problem \eqref{eq:lapp}--\eqref{eq:bep} provided $s$ is found from
the equation
\begin{equation}
\mathcal{R} (s) = r \, , \quad \mbox{where} \ \mathcal{R} (s) = [ s^2 - 2 \, \Omega
(1) + 2 \, h (s) ] / 3 \, . \label{eq:calR}
\end{equation}
The latter function has only one minimum, say $r_c > 0$, attained at some $s_c >
s_0$. Hence if $r \in (r_c, r_0)$, where
\[ r_0 = \lim_{s \to s_0 + 0} \mathcal{R} (s) = \frac{1}{3} \left[ s^2_0 -
2 \, \Omega (1) + 2 \, h_0 \right] ,
\]
then equation \eqref{eq:calR} has two solutions $s_+$ and $s_-$ such that $s_0 < s_+
< s_c < s_-$. By substituting $s_+$ and $s_-$ into \eqref{eq:Uim} and \eqref{eq:d},
one obtains the so-called {\it stream solutions} $(u_+, H_+)$ and $(u_-, H_-)$,
respectively. Indeed, these solutions satisfy Bernoulli's equation
\[ [ u'_\pm (H_\pm) ]^2 + 2 H_\pm = 3 r
\]
along with relations \eqref{eq:u}. It should be mentioned that $s_-$ and the
corresponding $H_-$ exist for all values of $r$ greater than $r_c$, whereas $s_+$
and $H_+$ exist only when $r$ is less than or equal to $r_0$; in the last case $s_+
= s_0$.

\subsection{Bounds for $(\psi, \Gamma)$} 

To express bounds for non-stream solutions of problem
\eqref{eq:lapp}--\eqref{eq:bep} we use solutions of problem \eqref{eq:u} and the
values $r_c$, $H_-$ and $H_+$; the last two serve as bounds for
\[ \hat{\Gamma} = \max_{(x, y) \in \Gamma} y \quad \mbox{and} \quad \check{\Gamma} =
\min_{(x, y) \in \Gamma} y .
\]
Now we formulate results generalising Theorems~1.1 and 1.2 in \cite{KK1} for
periodic solutions.

\begin{theorem} 
Let $(\psi, \Gamma)$ be a non-stream solution of problem
\eqref{eq:lapp}--\eqref{eq:bep} in the sense of Definition~1. Then the following
two assertions are true provided $\psi \leq 1$ on $\overline{D_\Gamma}$.

1. If $\check \Gamma < h_0$, then 
\begin{equation}
\psi (x, y) < U (y; \check s) \ in \ the \ strip \ \RR \times (0, \check
\Gamma) , \label{check}
\end{equation}
where $U$ is defined by formula \eqref{eq:Uim} and $\check s > s_0$ is such that $h
(\check s) = \check \Gamma$. Moreover, the inequalities $(A) \ r \geq r_c$, $(B) \
H_- \leq \check \Gamma$ hold, and if $r \leq r_0$, then $(C) \ \check \Gamma \leq
H_+$.

2. If $h_0 \neq +\infty$ and $\check \Gamma = h_0$, then inequality \eqref{check} is
nonstrict, whereas inequalities $(A)$--$(C)$ are true.
\end{theorem}

\begin{theorem}
Let $(\psi, \Gamma)$ be a non-stream solution of problem
\eqref{eq:lapp}--\eqref{eq:bep} in the sense of Definition~1. If $\hat \Gamma <
h_0$ and $\psi \geq 0$ on $\overline{D_\Gamma}$, then
\begin{equation}
\psi (x, y) > U (y; \hat s) \ in \ D_\Gamma , \label{hat}
\end{equation}
where $U$ is defined by formula \eqref{eq:Uim} and $\hat s > s_0$ is such that $h
(\hat s) = \hat \Gamma$. Moreover, $\hat \Gamma \geq H_+$ provided $r \leq r_0$ and
$\psi \leq 1$ on $\overline {D_\Gamma}$.
\end{theorem}

It occurs that some inequalities in Theorems~2 and 3 are strict under the assumption
that the first derivatives of $\psi$ are H\"older continuous near the points, where
the values $\check \Gamma$ and $\hat \Gamma$ are attained, as the following
assertions demonstrate.

\begin{proposition}
Let $(\psi, \Gamma)$ be a non-stream solution of problem
\eqref{eq:lapp}--\eqref{eq:bep} in the sense of Definition~1 and such that $\psi
\leq 1$ on $\overline{D_\Gamma}$. Also, let $\Gamma$ be of the class $C^{1,
\alpha}$, $\alpha \in (0, 1)$, near some point $(x^0, \check \Gamma) \in \Gamma$. If
$\psi \in C^{1, \alpha} (\bar X)$, where $X$ is the intersection of $D_\Gamma$ with
a sufficiently small neighbourhood of $(x^0, \check \Gamma)$, then the inequalities
are strict in $(A)$ and $(B)$, and if $r \leq r_0$, then the inequality in $(C)$ is
also strict.
\end{proposition}

\begin{proposition}
Let $(\psi, \Gamma)$ be a non-stream solution of problem
\eqref{eq:lapp}--\eqref{eq:bep} in the sense of Definition~1 and such that $\psi
\geq 0$ on $\overline{D_\Gamma}$. Also, let $\Gamma$ be of the class $C^{1,
\alpha}$, $\alpha \in (0, 1)$, near some point $(x^0, \hat \Gamma) \in \Gamma$. If
$\psi \in C^{1, \alpha} (\bar X)$, where $X$ is the intersection of $D_\Gamma$ with
a sufficiently small neighbourhood of $(x^0, \hat \Gamma)$, then $\hat \Gamma > H_+$
provided $r \leq r_0$ and $\psi \leq 1$ on $\overline {D_\Gamma}$.
\end{proposition}

\subsection{Proof of Theorem 2}

First, let $\check \Gamma < h_0$, and so there exists $\check s > s_0$ such that $h
(\check s) = \check \Gamma$, whereas the function $U (y; \check s)$ solves problem
\eqref{eq:u} on $(0, \check \Gamma)$. Moreover, formula \eqref{eq:Uim} defines $U
(y; \check s)$ on the half-axis $y \geq 0$ provided $\omega (t)$ is extended by $-1$
for $t > 1$. This implies that $\check s^2 > 2 \max_{\tau \geq 0} \Omega (\tau)$,
and we have
\[ U' (y; \check s) = \sqrt{\check s^2 - 2 \, \Omega (U (y; \check s))}
> 0 \quad \mbox{for} \ y \geq 0 .
\]
Hence $U (y; \check s)$ is a monotonically increasing function of $y$, $U (y; \check
s) > 1$ for $y > h_0$ and $U (y; \check s) \to +\infty$ as $y \to +\infty$.

Putting $U_\ell (y) = U (y + \ell; \check s)$ for $\ell \geq 0$, we see that $U_\ell
(y) > 1$ on $[0, \check \Gamma]$ when $\ell > \check \Gamma$. Therefore, $U_\ell -
\psi > 0$ on $\overline {D_\Gamma}$ for $\ell > \check \Gamma$. Let us show that there
is no $\ell_0 \in (0, \check \Gamma)$ such that
\begin{equation}
\min_{(x, y) \in \overline {D_\Gamma}} \{ U_{\ell_0} (y) - \psi (x, y) \} = 0 \, .
\label{inf}
\end{equation}

Assuming that such a value exists (in the case when there are several such values,
we denote by $\ell_0$ the largest of them), we see that \eqref{inf} holds only when
\[ U_{\ell_0} (y^0) - \psi (x^0, y^0) = 0 \quad \mbox{for some} \ (x^0, y^0) \in
\RR \times (0, \check \Gamma) 
\]
because $U_{\ell_0} - \psi$ is separated from zero on $\overline {D_\Gamma}
\setminus [\RR \times (0, \check \Gamma)]$. Moreover, \eqref{inf} implies that
\[ \nabla \psi = \nabla U_{\ell_0} \ \mbox{at} \ (x^0, y^0) , \ \mbox{and so} \
(x^0, y^0) \in D_\Gamma \setminus \Gamma_\psi .
\]
Since $D_\Gamma \setminus \Gamma_\psi$ is a domain, Theorem~1 is applicable in
$D_\Gamma \setminus \Gamma_\psi$, which yields that $\psi$ coincides with
$U_{\ell_0}$ there. Hence these functions coincide in $D_\Gamma$ because
$\Gamma_\psi$ has the zero measure. However, this contradicts to the fact that
$U_{\ell_0} - \psi$ is separated from zero on $\overline {D_\Gamma} \setminus [\RR
\times (0, \check \Gamma)]$.

The obtained contradiction shows that $U (y; \check s) \geq \psi (x,y)$ on
$\overline {D_\Gamma}$ and vanishes when $y = 0$. Moreover, Theorem~1 implies that
$U (\cdot; \check s)$ and $\psi$ cannot be equal at an inner point of $D_\Gamma
\setminus \Gamma_\psi$ because the latter function is not a stream solution.
Furthermore, $\nabla U \neq 0$ on $\Gamma_\psi$, and so the extended $U$ is strictly
greater than $\psi$ on $D_\Gamma$ which completes the proof of \eqref{check}.

To show that (A)--(C) are valid, we consider a point, say $(x^0, \check \Gamma)$, at
which the curve $\Gamma$ is tangent to $y = \check \Gamma$, and so $U (\check
\Gamma; \check s) - \psi (x^0, \check \Gamma) = 0$ because both terms on the
left-hand side are equal to one. It was proved that $U (y; \check s) - \psi (x, y)
\geq 0$ on $\RR \times [0, \check \Gamma]$, which implies
\begin{equation}
\big[ U' (y; \check s) - \psi_y (x,y) \big]_{(x,y)=(x^0,\check \Gamma)} \leq 0 \, .
\label{th_2}
\end{equation}
Since Bernoulli's equation at $(x_0, \check \Gamma)$ has the form $\psi_y (x^0,
\check \Gamma) = \sqrt{3 r - 2 \check \Gamma}$\,--\,cf. formula
\eqref{eq:ben}\,--\,inequality \eqref{th_2} gives
\[ U' (\check \Gamma; \check s) \leq \sqrt{3 r - 2 \check \Gamma} \quad 
\Longleftrightarrow \quad \check s^2 - 2 \Omega (1) \leq 3 r - 2 \check \Gamma .
\]
Hence $\mathcal{R} (\check s) \leq r$ in view of \eqref{eq:calR}, and combining the
latter inequality and $h (\check s) = \check \Gamma$, one obtains that (A) and (B)
are true in assertion~1. Moreover, (C) is also true provided $r \leq r_0$.

Now we turn to assertion 2 and begin with the case when $\check s = s_0 > 0$; here
the equality is a consequence of the assumption that $\check \Gamma = h_0$. Let us
introduce $\omega^{(\epsilon)} (\tau) = \omega (\tau) - \epsilon$, where $\epsilon >
0$ is small, and let $\Omega^{(\epsilon)} (\tau)$ and $U^{(\epsilon)} (y; s)$ be
defined by formulae \eqref{eq:Uim} with $\omega$ changed to $\omega^{(\epsilon)}$;
similarly, we define $h^{(\epsilon)} (s)$ using \eqref{eq:d}, whereas to obtain
$H^{(\epsilon)}_+$ and $H^{(\epsilon)}_-$ one has to combine \eqref{eq:calR} and
\eqref{eq:d}.

Since $s_0 > 0$, we have that
\[ s_0^{(\epsilon)} = \sqrt{\, 2 \max_{\tau \in [0, 1]} \Omega^{(\epsilon)} (\tau)}
< s_0 .
\]
Furthermore, it is straightforward to verify the inequalities
\[ h^{(\epsilon)} (s) < h (s) \quad \mbox{and} \quad U^{(\epsilon)} (y; s) > U (y; s)
\quad \mbox{for} \ s \geq s_0 .
\]
Therefore, $U^{(\epsilon)} (y; s_0)$ solves the problem on $(0, \check \Gamma)$
analogous to \eqref{eq:u}, but with $\omega$ is changed to $\omega^{(\epsilon)}$ and
with the value $U^{(\epsilon)} (\check \Gamma; s_0)$ greater than one. Moreover, in
view of the inequality
\[ U^{(\epsilon)''} + \omega \left( U^{(\epsilon)} \right) \geq 0 \quad \mbox{on} \
(0, \check \Gamma) ,
\]
Theorem 1 is applicable, and so the considerations used at the beginning of the
proof and involving $U (\cdot; \check s)$ are valid for $U^{(\epsilon)} (\cdot;
s_0)$ as well, thus yielding
\begin{equation}
\psi (x, y) < U^{(\epsilon)} (y; \check s) \quad \mbox{for} \ (x, y) \in \RR
\times (0, \check \Gamma) , \label{28}
\end{equation}
which is similar to \eqref{check}; here it is also taken into account that $\check
s = s_0$. Letting $\epsilon \to 0$ in the last inequality, one obtains
\[ \psi (x, y) \leq U (y; \check s) \quad \mbox{for} \ (x, y) \in \RR
\times (0, \check \Gamma) ,
\]
form which the inequalities in (A)--(C) follow in the same way as above.

Now let $s_0 = 0$. First we assume that $\Omega (\tau) < 0$ for $\tau \in (0, 1]$,
in which case we have
\[ U' (y; s_0) > 0 \ \mbox{for} \ y > 0 \quad \mbox{and} \quad U' (0; s_0) = 0 
\]
for the function defined by formulae \eqref{eq:Uim}. Then the considerations used in
the case when $\check \Gamma < h_0$ are applicable. Otherwise, the considerations
based on $U^{(\epsilon)} (y; s)$ yield \eqref{28}, and the results follow letting
$\epsilon \to 0$.

\subsection{Proof of Theorem 3} 

At its initial stage the proof of this theorem is similar to that of Theorem~2.
Namely, we consider the case when $\hat \Gamma < h_0$ first. Since there exists
$\hat s > s_0$ such that $h (\hat s) = \hat \Gamma$, the function $U$ given by
formula \eqref{eq:Uim} with $s = \hat s$ solves problem \eqref{eq:u} on $(0, \hat
\Gamma)$. Moreover, the same formula defines this function for all $y \leq \hat
\Gamma$ provided $\omega (t)$ is extended to $t < 0$ by zero, in which case $\hat
s^2 > 2 \max \Omega (\tau)$. Then
\[ U' (y; \hat s) = \sqrt{\hat s^2 - 2 \, \Omega (U (y; \hat s))} > 0 \quad 
\mbox{for} \ y \leq \hat \Gamma ,
\]
and so $U (y; \hat s)$ is a monotonically increasing function of $y$ such that $U
(y; \hat s) < 0$ for $y < 0$.

Let $U_\ell (y) = U (y - \ell; \hat s)$ for $\ell \geq 0$, which implies that
$U_\ell (y) < 0$ on $[0, \hat \Gamma]$ provided $\ell > \hat \Gamma$. Therefore,
$U_\ell - \psi < 0$ on $\overline {D_\Gamma}$ for $\ell > \hat \Gamma$. Similarly to
the proof of Theorem~2, one obtains that there is no $\ell_0 \in (0, \hat \Gamma)$
such that
\[ \max_{(x, y) \in \RR \times [0, \hat \Gamma]} \{ U_{\ell_0} (y) - \psi (x,
y) \} = 0 \, . 
\]
Hence $U (y; \hat s) - \psi (x,y)$ is non-positive on $\overline {D_\Gamma}$ and
vanishes when $y = 0$. Now, applying Theorem~1 in the same way as in the proof of
Theorem~2, we arrive at inequality \eqref{hat}.

To prove that $\hat \Gamma \geq H_+$, we argue by analogy with the proof of
Theorem~2. In view of periodicity of $\Gamma$, there exists $(x^0, y^0) \in \Gamma$
such that $y^0 = \hat \Gamma$ (it is clear that $\Gamma$ is tangent to $y = \hat
\Gamma$ at this point). Then $U (\hat \Gamma; \hat s) - \psi (x^0, \hat \Gamma) = 0$
because both terms on the left-hand side are equal to one. Since $\psi$ is a
non-stream solution, then
\[ \big[ U' (y; \hat s) - \psi_y (x,y) \big]_{(x,y)=(x^0,\hat \Gamma)} \geq 0 \, .
\]
Using Bernoulli's equation for $\psi$ at $(x^0, \hat \Gamma)$, we show that
\begin{equation}
 U' (\hat \Gamma; \hat s) \geq \sqrt{3 r - 2 \hat \Gamma} . \label{20}
\end{equation}
Indeed, it follows from the boundary condition $\psi (x^0, \hat \Gamma) = 1$ and the
assumption that $\psi \leq 1$ on $\overline{D_\Gamma}$ that $\psi_y (x^0, \hat
\Gamma)$ is non-negative, and so
\[ \psi_y (x^0, \hat \Gamma) = \sqrt{3 r - 2 \hat \Gamma} .
\]
Combining this and the inequality preceding \eqref{20}, we see that \eqref{20} is
true, which implies the required inequality.

\subsection{Proof of Propositions 1 and 2}

Since the proof of Proposition 2 is similar to that of Proposition 1, we restrict
ourselves to proving the latter assertion only.

To prove Proposition~1, we notice that there exists $\check s > s_0$ such that $h
(\check s) = \check \Gamma$ and $U' (y; \check s) > 0$ for $y \geq 0$; here $U
(\cdot; \check s)$ is defined by formula \eqref{eq:Uim} provided $\omega (t)$ is
extended to $t > 1$ by zero. (This follows in the same way as in the proof of
Theorem~2.) Since inequality \eqref{th_2} is valid, it yields that $\psi_y (x^0,
\check \Gamma) > 0$.

Now we apply the partial hodograph transform in two neighbourhoods of $(x^0, \check
\Gamma)$ so that in both cases the image of this point is $(q^0, 1)$ on the $(q,
p)$-plane. In the first case, some $X_1 \subset D_\Gamma$ is mapped to a
neighbourhood 
\[ Q_1 \subset \{ (q, p): q \in \RR , \, p < 1 \} ,
\]
and $h_1 (q, p)$ in $Q_1$ corresponds to $\psi (x, y)$ defined in $X_1$. In the
second case, some $X_2 \subset \RR \times (0, 1)$ is mapped to another
neighbourhood $Q_2$ on the same plane as $Q_1$, whereas $h_2 (q, p)$ in $Q_2$
corresponds to $U (y; \check s)$. Since inequality \eqref{check} holds for $\psi$
and $U (\cdot; \check s)$, we have that $h_1 - h_2 > 0$ in some subset of $Q_1 \cap
Q_2$ which is the intersection of $\{ (q, p): q \in \RR , \, p < 1 \}$ with a
neighbourhood of $(q^0, 1)$. Besides, $h_1 - h_2$ vanishes at $(q^0, 1)$ because
$\psi (x^0, \check \Gamma) = U (\check \Gamma; \check s) = 1$. Then it follows from
Lemma 2 that
\[ [ \partial_p (h_1 - h_2) ]_{(q, p)=(q^0, 1)}  < 0 ,
\]
which implies that inequality \eqref{th_2} is strict. Using this fact in the
considerations that follow \eqref{th_2}, we obtain that the inequalities in
(A)--(C) are strict.

\subsection{Discussion}

In his renown book \cite{C}, the first edition of which was published in 1918,
Carath\'eodory had proved a quite general theorem for the first order ordinary
differential equation. It concerns the existence of a solution which satisfies the
equation on an interval up to a set of Lebesgue-measure zero. The proof is based on
assumptions whose general form is now referred to as the Carath\'eodory condition
(see \cite{AZ}, \S\,1.4, for its discussion). In the framework of this condition,
$f$ is supposed to be continuous in almost all papers, dealing with equation
\eqref{eq}, inequality \eqref{in} and their generalisations (see, for example, the
notes \cite{K} and \cite{O} by Keller and Osserman, respectively, dating back to
1957, and numerous papers citing these notes). It is also worth mentioning in this
connection, that non-uniqueness takes place for the first order ordinary
differential equation when the smoothness of a nonlinear term is less than Lipschitz
with respect to the unknown function.

Furthermore, considering equation \eqref{eq} in \cite{GNN} (see Theorem~1 on
p.~209), the authors require even more, namely, that $f$ is of class $C^1$. This
substantially simplifies treatment of the equation comparing with Theorem~1 in the
present paper, where $f \in L^p_{loc} (\RR)$ for $p > n$. On the other hand, the
assumption imposed on solutions in our theorem, namely, that their gradients are
non-vanishing, is essential. This condition allows us to avoid non-uniqueness even
without the Carath\'eodory condition.

Turning to the problem of periodic water waves with vorticity, the papers \cite{CS1}
and \cite{MM} should be mentioned. Discontinuous vorticity distributions from
$L^\infty$ are considered in the first of them, whereas the distribution is merely
$L^p$-integrable with an arbitrary $p \in (1, \infty)$ in the second one. However,
only unidirectional flows (they have no stagnation points within the fluid) are
studied in both papers, and in this case the global partial hodograph transform can
be applied to simplify the problem, thus reducing the effect of non-smooth
vorticity.

It is worth mentioning that the assumptions on a solution of the water wave problem
here are weaker not only than those imposed in our recent paper \cite{KK1}, but also
than those in \cite{KKL1}. Indeed, the most restrictive condition in \cite{KKL1} is
that the horizontal component of the velocity field is bounded from below by a
positive constant.

\vspace{6mm}

\noindent {\bf Acknowledgements.} V.~K. was supported by the Swedish Research
Council (VR). N.~K. acknowledges the support from the G.~S.~Magnuson's Foundation of
the Royal Swedish Academy of Sciences and Link\"oping University.

\end{document}